\documentclass{amsart}
\usepackage{graphicx}
\usepackage{amscd}

\newtheorem{theorem}{Theorem}
\theoremstyle{plain}

\newtheorem{corollary}{Corollary}

\newtheorem{definition}{Definition}
\newtheorem{example}{Example}

\newtheorem{lemma}{Lemma}

\newtheorem{remark}{Remark}

\numberwithin{equation}{section}

\begin{document}
\title[Compactly almost periodic homeomorphisms]{Characterizations of compactly almost periodic homeomorphisms of metrizable
spaces}
\author{Paul Fabel}
\address{Drawer MA Mississippi State, MS 39762}
\email{fabel@ra.msstate.edu}
\urladdr{http://www.msstate2.edu/\symbol{126}fabel}
\date{August 28 2001}

\begin{abstract}
For metrizable spaces we replace the notion of almost periodic homeomorphism
with a similar notion and verify that the usual characterizations of almost
periodic homeomorphisms of compact metric spaces are valid for all
metrizable spaces. We include examples and prove some related results.
\end{abstract}

\maketitle

\section{\protect\bigskip Introduction}

A homeomorphism $h:X\rightarrow X$ of a metric space $(X,d)$ is \textbf{%
almost periodic} if given $\varepsilon >0$ there exists an integer $N$ so
that each block of $N$ consecutive iterates of $h$ contains a map $h^{n}$
such that $\forall x\in X$ $d(x,h^{n}(x))<\varepsilon .$

Classifying such maps for various choices of $(X,d)$ is the subject of \cite
{Brech},\cite{Foland},\cite{lex}, and \cite{Ritter}. A useful tool is the
result of Gottschalk \cite{Gott}, \cite{Sam}.

\begin{theorem}[Gottschalk]
\label{Gotts}If $X$ is a compact metric space then a homeomorphism $%
h:X\rightarrow X$ is almost periodic if and only if $\{h^{n}|n\in Z\}$ is an
equicontinuous family of maps.
\end{theorem}

Unfortunately the theorem is false if $X$ is not compact since the property
of almost periodicity is metric dependent. For example if $h$ is an
irrational rotation of the plane about the origin then $h$ is not almost
periodic with the usual metric, but $h$ is almost periodic when $R^{2}$ is
equipped with a metric inherited from the compact $2$ sphere $R^{2}\cup
\{\infty \}$.( See Theorem \ref{apv2}, example \ref{ap3}, \cite{Brech} and 
\cite{lex}).

The main result of this paper, Theorem \ref{mainthm}, serves to resurrect
Theorem \ref{Gotts} by replacing the notion of almost periodic with the
notion of a \textbf{compactly almost periodic }homeomorphism defined as
follows. If $X$ is a metrizable space then a homeomorphism $h:X\rightarrow X$
is \textbf{compactly almost periodic }if the orbit closure of each compactum 
$B\subset X$ under the action of $h$ is compact, and if for each compact
invariant subset $A\subset X,$ $h_{A}:A\rightarrow A$ is almost periodic.

Though generally distinct, the two notions are often comparable:

\begin{enumerate}
\item  Almost periodic homeomorphisms of locally complete metric spaces are
compactly almost periodic.

\item  If $X$ is compact then $h$ is compactly almost periodic if and only
if $h$ is almost periodic for each metric on $X.$

\item  If $X$ is separable and locally compact then $h$ is compactly almost
periodic if and only if $h$ is almost periodic for some metric on $X.$
\end{enumerate}

In particular a homeomorphism $h:R^{2}\rightarrow R^{2}$ is compactly almost
periodic if and only if $h$ is conjugate to rotation or reflection (
Corollary \ref{planecap})

As in the compact case, $h$ is compactly almost periodic if and only if
there exists a metric $d^{\ast }$ on $X$ such that $\overline{\{h^{n}\}}$ is
a compact abelian group of isometries of $X$ ( the closure is taken in $%
C(X,X)$ with the compact open topology). Consequently, the space of orbit
closures under $h$ forms an uppersemicontinuous decomposition of $X$ (
compatible with the Hausdorf metric) into compacta each of which is a
topological abelian group ( Theorem \ref{Decomp}).

One difficulty created by the failure of $X$ to be compact includes the
following: There exists a compactly almost periodic isometry of a locally
compact metric space $X$ such that $\overline{\{h^{n}\}}$ is a compact
isometry group but $\overline{\{h^{n}\}}$ is not metrizable. (See example 
\ref{nometric}).

All function spaces are endowed with the compact open topology and all
references to metrics $d$ on metrizable spaces $X$ are assumed (or proven)
to induce a topology compatible with that of $X.$

\section{Main Result}

\begin{theorem}
\label{mainthm} Suppose $X$ is a metrizable space and $h:X\rightarrow X$ is
a homeomorphism such that $\overline{\{h^{n}(x)\}}$ is compact $\forall x\in
X.$ Then the following are equivalent:
\end{theorem}

\begin{enumerate}
\item  $h$ is a compactly almost periodic homeomorphism.

\item  For some metric $d$ on $X$ $\{h^{n}\}$ is an equicontinuous family of
maps.

\item  For some metric $d^{\ast }$ on $X$ $h:(X,d^{\ast })\rightarrow
(X,d^{\ast })$ is an isometry.

\item  For each metric $d$ on $X$ $\{h^{n}\}$ is an equicontinuous family of
maps.

\item  $\overline{\{h^{n}\}}$ is a compact subspace of $C(X,X)$.

\item  There exists a metric $d^{\ast }$ on $X$ such that $\overline{%
\{h^{n}\}}\subset C(X,X)$ is a compact abelian topological group of
isometries of $X.$
\end{enumerate}

\begin{proof}
Trivially $4\Rightarrow 2.$ $2\Rightarrow 5$ and $5\Rightarrow 4$ by Theorem 
\ref{Ascoli1}. Thus $2,4,$ and $5$ are equivalent. By Lemma \ref{equiiso} $%
2\Rightarrow 3.$ Trivially $3\Rightarrow 2$. $1$ and $4$ are equivalent by
Lemma \ref{apiffeq}. Finally, $6\Rightarrow 5$ is trivial and $5\Rightarrow
6 $ by Theorem \ref{Bcmptgrp}.
\end{proof}

\section{Definitions}

\begin{definition}
If $A$ is a metric space then $2^{A},$ the collection of all compact subsets
of $A$ forms a metric space with the \textbf{Hausdorf} metric generated by
the following condition. Two compacta $\{B,C\}\subset 2^{A}$are less than a
distance $\varepsilon $ from each other if each point of $b\in B$ is within $%
\varepsilon $ of some point of $C$ and if each point $c\in C$ is within $%
\varepsilon $ of some point of $B.$
\end{definition}

\begin{definition}
Let $C(A,Y)$ denote the space of continuous function from the metric space $%
A $ to the metric space $Y$ with the compact open topology. Let $%
H(X,X)\subset C(X,X)$ denote the space of homeomorphisms from $X$ to $X.$
\end{definition}

\begin{definition}
If $(Y,d_{Y})$ is a metric space then a basis for $C(X,Y)$ consists of sets $%
B_{A}(f,\varepsilon )$ where $A\subset X$ is compact and $g\in
B_{A}(f,\varepsilon )$ iff $d_{Y}(f(a),g(a))<\varepsilon $ $\forall a\in A.$
For $\{f,g\}\subset C(X,Y)$ let $d_{Y}(f,g)=\sup_{x\in X}d(f(x),g(x)).$ Note 
$d_{Y}(f,g)\leq \infty .$
\end{definition}

\begin{definition}
Suppose $(X,d)$ is a metric space. Then a homeomorphism $h:X\rightarrow X$
is \textbf{almost periodic} if for each $\varepsilon >0$ there exists $N\geq
0$ so that if $M\in Z$ then there exists $i$ such that $M\leq i\leq N+M$ and
such that $\forall x\in X$ $d(h^{i}(x),x)<\varepsilon .$ If $X$ is a
metrizable space and $h:X\rightarrow X$ is a homeomorphism such that $%
\overline{\cup _{n\in Z}h^{n}(B)}$ is compact for each compact $B\subset X,$
and if for each compact invariant set $A,$ $h_{A}:A\rightarrow A$ is an
almost periodic homeomorphism of $A$ then $h$ is \textbf{compactly almost
periodic. }
\end{definition}

\begin{definition}
Suppose $(A,d_{A})$ and $(Y,d_{Y})$ are metric spaces. A subset $\mathcal{F}$
$\subset C(A,Y)$ is \textbf{equicontinuous} if for each $a\in A$ and each $%
\varepsilon >0$ there exists $\delta >0$ such that for each $f\in \mathcal{F}
$, if $d_{A}(a,x)<\delta $ then $d_{Y}(f(a),f(x))<\varepsilon .$
\end{definition}

\begin{definition}
If $\mathcal{F}\subset C(X,Y)$ and $A\subset X$ then $\mathcal{F}_{A}\subset
C(A,Y)$ denotes all maps $f_{A}$ such that $f\in \mathcal{F}.$
\end{definition}

\begin{definition}
If $(X,d)$ is a metric space then a homeomorphism $h:X\rightarrow X$ is an 
\textbf{isometry }if $d(x,y)=d(h(x),h(y)).$ Let $ISOM(X,d)\subset C(X,X)$
denote the collection of isometries from $X$ to $X.$
\end{definition}

\begin{definition}
The metric space $(X,d)$ is \textbf{locally complete} if for each $x\in X$
there exists $\varepsilon >0$ such that $\overline{B(x,\varepsilon )}$ is
complete with respect to $d.$
\end{definition}

\section{Comparing almost periodic and compactly almost periodic
homeomorphisms}

The notions `almost periodic homeomorphism' and `compactly almost periodic
homeomorphism' are distinct but often comparable.

\subsection{Positive comparisons}

\begin{lemma}
\label{isocap}Suppose $(X,d)$ is a locally complete metric space and $%
h:X\rightarrow X$ is an almost periodic isometry. Then $h$ is compactly
almost periodic.
\end{lemma}

\begin{proof}
Since $h$ is an isometry $\{h^{n}\}$ is equicontinuous. By Theorem \ref
{mainthm} it suffices to prove $\overline{\{h^{n}(x)\}}$ is compact $\forall
x\in X.$ Suppose $x\in X.$ Choose $\gamma >0$ so that $\overline{B(x,\gamma )%
}$ is complete. To prove $\overline{\{h^{n}(x)\}}$ is complete it suffices
to prove that each Cauchy sequence of the form $\{h^{n_{k}}(x)\}$ has a
limit in $X$. Suppose $\{h^{n_{k}}(x)\}$ is Cauchy. Since $\{h^{n_{k}}(x)\}$%
is Cauchy we may choose $K$ so that $d(h^{n_{m}}(x),h^{n_{l}}(x))<\gamma $
if $m,l\geq K.$ Since $h^{-N_{K}}$ is an isometry $\{h^{n_{k}-N_{K}}(x)\}$
is also Cauchy and

$d(h^{n_{m}-N_{K}}(x),h^{n_{l}-N_{K}}(x))<\gamma $ if $m,l\geq K.$ Taking $%
m=K$ $d(h^{n_{m}-N_{K}}(x),h^{n_{l}-N_{K}}(x))=d(x,h^{n_{l}-N_{K}})<\gamma .$
Consequently $\{h^{n_{k}-N_{K}}(x)\}$ has a limit $y$. Hence $%
\{h^{n_{k}}(x)\}$ has limit $h^{N_{K}}(y).$ Thus $\overline{\{h^{n}(x)\}}$
is complete. By Lemma \ref{tb} $\overline{\{h^{n}(x)\}}$ is totally bounded
and therefore $\overline{\{h^{n}(x)\}}$ is compact.
\end{proof}

\begin{theorem}
\label{apv1}Suppose $(X,d)$ is a metric space and $h:X\rightarrow X$ is an
almost periodic homeomorphism. Then $h$ is compactly almost periodic if and
only if $\forall x\in X$ $\overline{\{h^{n}(x)\}}$ is locally complete with
respect to $d.$
\end{theorem}

\begin{proof}
If $h$ is compactly almost periodic then $\overline{\{h^{n}(x)\}}$ is
compact and hence complete and therefore locally complete. Conversely
suppose $\overline{\{h^{n}(x)\}}$ is locally complete $\forall x\in X.$ By
Lemma \ref{apequi} $\{h^{n}\}$ is equicontinuous. By Lemma \ref{tb} $%
\overline{\{h^{n}(x)\}}$ is totally bounded and hence bounded. Let $d^{\ast
}(x,y)=\sup_{n}d(h^{n}(x),h^{n}(y)).$ By Lemma \ref{equiiso} $d^{\ast }$ is
a metric equivalent to $d$ such that $h$ is an isometry of $(X,d^{\ast }).$
Suppose $x\in X.$ Let $Y=\overline{\{h^{n}(x)\}}.$ To prove $Y$ is locally $%
d^{\ast }$ complete suppose $y\in Y.$ Choose $\gamma >0$ so that $A=%
\overline{B_{d}(y,\gamma )}\cap \overline{\{h^{n}(x)\}}$ is $d$ complete.
Choose $\beta >0$ so that $d(y,z)<\beta \Rightarrow d^{\ast }(y,z)<\gamma .$
Let $C=\overline{B_{d^{\ast }}(y,\beta )}\cap \overline{\{h^{n}(x)\}}.$ Note 
$C\subset A$ and $C$ and $A$ are closed. If $\{z_{n}\}$ is $d^{\ast }$
Cauchy in $C$ then $\{z_{n}\}$ is $d$ Cauchy in $A.$ Thus $z=\lim z_{n}$
exists and $z\in A\cap C.$ Hence $C$ is $d^{\ast }$ complete and therefore $%
Y $ is locally $d^{\ast }$ complete. Now we will prove that $h_{Y}$ is an
almost periodic isometry of $Y.$ Suppose $\varepsilon >0.$ Choose $\delta
<\varepsilon $ so that $d(x,y)<\delta \Rightarrow
d(h^{n}(x),h^{n}(y))<\varepsilon /2.$ Choose $N$ so that each block of
iterates contains a map $h^{n}$ such that $d(h^{n},id)<\delta .$ Suppose $%
n\in Z.$ Choose $0\leq m\leq N$ such that $d(h^{n-m},id)<\delta .$ Thus $%
d(h^{n-m}(x),x)<\delta .$ Hence $d^{\ast }(h^{n-m}(x),x)<\varepsilon /2.$
However $h^{n-m}$ is a $d^{\ast }$ isometry. Thus $d^{\ast
}(h^{n-m}(h^{k}(x)),h^{k}(x))<\varepsilon /2.$ Hence, since $\{h^{k}(x)\}$
is dense in $Y,$ $d^{\ast }(h_{Y}^{n-m},id_{Y})\leq \varepsilon
/2<\varepsilon .$ By Lemma \ref{invariant} $h(Y)=Y.$ Therefore $h_{Y}$ is an
almost periodic isometry of the locally complete space $Y.$ Thus by Lemma 
\ref{isocap} $h_{Y}$ is compactly almost periodic. In particular, since $%
\{h_{Y}^{n}(x)\}$ is dense in $Y,$ $Y$ is compact. Hence by Theorem \ref
{Ascoli1} $\overline{\{h^{n}\}}$ is compact. Therefore by Theorem \ref
{mainthm} $h$ is compactly almost periodic.
\end{proof}

\begin{corollary}
\label{apvcap}If $(X,d)$ is a locally compact metric space then each almost
periodic homeomorphism of $X$ is compactly almost periodic.
\end{corollary}

\begin{theorem}
\label{apv2}Suppose $X$ is a locally compact separable metrizable space and $%
h:X\rightarrow X$ is a homeomorphism. Then $h$ is compactly almost periodic
if and only if there exists a metric $d$ on $X$ such that $h$ is an almost
periodic homeomorphism of $(X,d).$
\end{theorem}

\begin{proof}
Suppose $h$ is compactly almost periodic. By Theorem \ref{mainthm} choose a
metric $d^{\ast }$ on $X$ such that $H=\overline{\{h^{n}\}}\subset C(X,X)$
is a compact isometry group of $X.$ Let $Y=X\cup \{\infty \},$ the one point
compactification of $X$ with metric $d.$ By Remark \ref{metrizable} $Y$ is
metrizable. By Theorem \ref{Bcmptgrp} and Remark \ref{hyygrp} $H$ and $%
H(Y,Y) $ are topological groups. Define $\phi :H\hookrightarrow H(Y,Y)$ as $%
\phi (g)=G$ where $G:Y\rightarrow Y$ is the unique homeomorphism such that $%
G_{|X}=g$ and $G(\infty )=\infty .$ Continuity of $\phi $ will follow from
continuity of $\phi $ at $id_{X}.$ Let $W=B_{Y}(id_{Y},\varepsilon )\subset
H(Y,Y).$ Let $A=Y\backslash B_{d^{\ast }}(\infty ,\varepsilon /3).$ Let $%
C=Y\backslash B_{d^{\ast }}(\infty ,2\varepsilon /3).$ Let $U=\{p\in
G|p(C)\subset int(A)\}.$ Let $V_{1}$ be an open set containing $id_{X}$ such
that $V_{1}^{-1}\subset U.$ Let $V_{2}=\{g\in G|d^{\ast
}(id_{A},g_{A})<\varepsilon \}.$Let $V=V_{1}\cap V_{2}.$ Suppose $g\in V.$
If $x\in A$ then $d^{\ast }(x,g(x))<\varepsilon $ since $g\in V_{2}.$ If $%
x\notin A$ then $g(x)\notin C$ since $g^{-1}\in U.$ Thus $d^{\ast
}(g(x),x))\leq d^{\ast }(x,\infty )+d^{\ast }(\infty ,g(x))<\varepsilon
/3+2\varepsilon /3=\varepsilon .$ Hence $g(V)\subset U.$ Therefore $\phi $
is continuous. Hence $\phi (H)=\overline{\{\phi (h)^{n}\}}$ is compact and
therefore by Theorem \ref{mainthm} $\phi (h)$ is an almost periodic
homeomorphism of $(Y,d).$ Thus, since $\phi (h)(X)=X$ and $\phi
(h)_{|X}=h_{|X}$ it follows that $h$ is an almost periodic homeomorphism of $%
(X,d).$ The converse follows from Corollary \ref{apvcap}.
\end{proof}

\begin{corollary}
\label{planecap}Orientation preserving compactly almost periodic
homeomorphisms of the plane $R^{2}$ are conjugate to rotation about the
origin. Orientation reversing compactly almost periodic homeomorphisms of $%
R^{2}$ are conjugate to reflection about the $y$ axis.
\end{corollary}

\begin{proof}
If $h:R^{2}\rightarrow R^{2}$ is compactly almost periodic then, following
the proof of Theorem \ref{apv2}, $h$ extends to an almost periodic
homeomorphism of $S^{2}=R^{2}\cup \{\infty \}$ fixing $\{\infty \}.$ Now
apply the results of Ritter \cite{Ritter}.
\end{proof}

\subsection{Negative comparisons and other counterexamples}

In each of the following examples $X$ is a metrizable space and $%
h:X\rightarrow X$ is either a homeomorphism or an isometry. Examples \ref
{nometric} and \ref{noclosed} illustrate two `pathological' isometry groups.
In examples \ref{ap1} and \ref{ap2} $\ h$ fails to be compactly almost
periodic despite some favorable properties.

\begin{example}
Let $S^{1}\subset R^{2}$ be the unit circle with the usual metric. Let $%
h:S^{1}\rightarrow S^{1}$ be an irrational rotation. Let $X=\{\cup _{n\in
Z}h^{n}(1)\}\subset S^{1}.$ Then $h_{|X}:X\rightarrow X$ is an almost
periodic homeomorphism but $h$ is not compactly almost periodic.
\end{example}

\begin{example}
\label{nometric}The homeomorphism $h:X\rightarrow X$ is compactly almost
periodic, but the compact isometry group $\overline{\{h^{n}\}}$ is not
metrizable, and $h$ fails to be almost periodic for each metric on $X.$
\end{example}

\begin{proof}
Let $J\subset S^{1}$ be an uncountable set such that for each finite subset $%
A=\{\alpha _{1},..\alpha _{N}\}\subset J,$ each point $\{z_{1},...z_{n}\}\in
\Pi _{i=1}^{N}S^{1}$ under the homeomorphism $r_{A}:\Pi
_{i=1}^{N}S^{1}\rightarrow \Pi _{i=1}^{N}S^{1}$ defined as $%
r_{A}(x_{1},...,x_{N})=(\alpha _{1}x_{1},..,\alpha _{N}x_{N}).$ Let $X$
denote the disjoint union of uncountably many copies of $S^{1}$indexed by $J$%
. Let points of distinct circles in $X$ have distance $1$ and let each
circle $S_{j}^{1}$be isometric to $S^{1}.$ For $j\in J$ define $%
h_{j}:S_{j}^{1}\rightarrow S_{j}^{1}$ to be the rotation $h_{j}(x)=\alpha
_{j}x.$ Define $h:X\rightarrow X$ as $h=\cup _{j\in J}h_{j}.$ Note $h$ is a
compactly periodic isometry of $X$. Hence $\overline{\{h^{n}\}}$ is a
compact group of isometries of $X$ and $g(S_{j})=S_{j}\forall j\in J$ $%
\forall g\in \overline{\{h^{n}\}}.$ Let $Y=\Pi _{j\in J}S_{j}^{1}.$ Choose a
`basepoint' $x_{j}$ from each $S_{j}^{1}\subset X.$ Define $\phi :\overline{%
\{h^{n}\}}\rightarrow Y$ as $\phi (g)=\{g(x_{j})\}.$ Then $\phi $ maps onto
a dense subset of $Y$ and hence $\phi $ is surjective. Thus $\overline{%
\{h^{n}\}}$ is not metrizable. Consequently for each equivalent metric $d$
on $X,$ $h$ is not an almost periodic homeomorphism of $(X,d)$ since if $h$
were an almost periodic homeomorphism of $(X,d)$ then $\overline{\{h^{n}\}}$
would be metrizable via the uniform metric. See \cite{Fabel} for more
details.
\end{proof}

\begin{example}
\label{noclosed}Isometry groups need not be closed in $C(X,X).$
\end{example}

\begin{proof}
Let $X=\{1,2,3,..\}$ with metric $d(m,n)=1$ iff $m\neq n.$ Let $%
f_{n}:X\rightarrow X$ be a bijection such that $f_{n}(k)=2k$ if $k\leq n.$
Define $f:X\rightarrow X$ as $f(n)=2n.$ Then $f_{n}\rightarrow f$, $f_{n}$
is an isometry of $X$ but $f$ is not an isometry of $X.$ Hence $ISOM(X,d)$
is not a closed subspace of $C(X,X).$
\end{proof}

\begin{example}
\label{ap1}$h:X\rightarrow X$ is a homeomorphism, $\overline{h^{n}(x)}$ is
compact $\forall x\in X$ and $h_{A}:A\rightarrow A$ is almost periodic for
each compact invariant set $A\subset X,$but $h$ is not compactly almost
periodic.
\end{example}

\begin{proof}
Let $\theta \in S^{1}.$ Let $X=((\{...1/3,1/2,1\}\times S^{1})\cup
\{0,\theta \}.$ Let $h(1/n,z)=(1/n,e^{\frac{i}{n}}z)$ for $n\in Z^{+}$ and $%
h(0,\theta )=(0,\theta ).$
\end{proof}

\begin{example}
\label{ap2}$X$ is compact,$\forall x\in X$ if $A=\overline{h^{n}(x)}$ then $%
h_{A}:A\rightarrow A$ is an almost periodic isometry, but $h$ is not
compactly almost periodic.
\end{example}

\begin{proof}
Let $X=D^{2},$ the closed unit disk in $R^{2}.$ Let $h(r,\theta
)=(r,e^{ir}\theta ).$
\end{proof}

\begin{example}
\label{ap3}Let $X=R^{2}.$ Define $h:R^{2}\rightarrow R^{2}$ as $h(z)=e^{2\pi
i\theta }z$ where $\theta \in R\backslash Q.$ Then $h$ is compactly almost
periodic but $h$ is not almost periodic with the usual metric on $R^{2}$. If 
$d$ is any metric on $R^{2}\cup \{\infty \}$ then $h$ is almost periodic
with respect to $d.$ See Theorem \ref{apv2}.
\end{example}

\section{Application: A continuous decomposition of $X$ into compact
metrizable abelian groups}

A decomposition $X^{\ast }$ of a topological space $X$ is a collection of
disjoint closed subsets of $X$ whose union is $X.$ Endowing $X^{\ast }$ with
the quotient topology, the natural quotient map $\pi :X\rightarrow X^{\ast }$
is both open and closed if and only if $X^{\ast }$ is both \textbf{%
uppersemicontinuous and lowersemicontinuous}\cite{Daver}\textbf{.} Such a
decomposition is called \textbf{continuous. }

\begin{theorem}
\label{Decomp}Suppose $X$ is a metrizable space, $h:X\rightarrow X$ is a
compactly almost periodic homeomorphism. Then sets of the form $\overline{%
\{h^{n}(x)\}}$ determine a continuous decomposition of $X$ into compact
abelian groups.
\end{theorem}

\begin{proof}
Let $H=\{h^{n}\}$ and let $G=\overline{H}\subset C(X,X)$. By theorem \ref
{Bcmptgrp} $G$ is a compact abelian topological group of isometries of $X$
for some metric $d.$ By Lemma \ref{orbitequals} $G(x)=\overline{H(x)}.$ Thus
if $\{x,y\}\subset X$ then $y\in \overline{H(x)}$ iff there exists $g\in G$
such that $y=g(x)$ iff $G(x)=G(y).$ Thus the collection of orbit closures
under the full action of $h$ forms a partition of $X$ into disjoint
compacta. Suppose $z^{\ast }\in G(x)$ and $w^{\ast }\in G(y)$. Then for each 
$g\in G$ $d(z^{\ast },w^{\ast })=d(g(z^{\ast }),g(w^{\ast })).$ Thus the
Hausdorf distance between $G(x)$ and $G(y)$ is the minimum of $d_{z\in
G(x),w\in G(y)}d(z,w).$ Hence, endowing $X/G$ with the Hausdorf topology,
the natural map $\pi :X\rightarrow X/G$ is both continuous and open since $%
d(x,y)\geq d(\pi (x),\pi (y))).$ In particular $\pi $ is a quotient map. To
show $\pi $ is also a closed map suppose $F\subset X$ is closed and $A=\cup
_{G(x)\cap F\neq \emptyset }G(x).$ Suppose $a\in \overline{A}$ choose $%
a_{n}\in A$ such that $a_{n}\rightarrow a.$ Let $a_{n}=g_{n}(x_{n})$ where $%
g_{n}\in G$ and $x_{n}\in F.$ Since $G$ is compact let $g_{n_{k}}\rightarrow
g\in G.$ By Lemma \ref{closedmult} $g_{n_{k}}^{-1}\rightarrow g^{-1}\in G.$
Let $B=a\cup \{a_{n_{k}}\}.$ $B$ is compact. Thus $g_{n_{k}}^{-1}(a_{n_{k}})%
\rightarrow g^{-1}(a).$ Thus $g^{-1}(a)\in F$ since $F$ is closed. Let $%
x=g^{-1}(a).$ Hence $a=g(x).$ Thus $\overline{A}=A$. Hence $\pi $ is a
closed map. Hence the orbit closures of $H$ determine continuous
decomposition of $X.$ Now we check the group structure of each orbit $G(x).$
Let $K(x)=\{g\in G|g(x)=x\}.$ Note $K(x)$ is a compact abelian ( normal)
subgroup of $G.$ Hence the compact group $G/K(x)$ acts freely and
effectively on $G(x).$ Hence $G(x)$ is a compact topological abelian group
under the multiplication $g_{1}(x)\ast g_{2}(x)=g_{1}g_{2}(x)$.
\end{proof}

\section{Theorems Lemmas and Remarks}

\subsection{Arzela-Ascoli Theorem}

\begin{theorem}[Arzela-Ascoli]
\label{Ascoli1} Suppose each of $X$ and $Y$ is a metric space and $\mathcal{%
F\subset }C(X,Y)$ in the compact open topology$.$ Suppose $\forall x\in X$ $%
\{f(x)|f\in \mathcal{F}\}$ has compact closure in $Y$. The following are
equivalent
\end{theorem}

\begin{enumerate}
\item  $\mathcal{F}$ is equicontinuous

\item  $\mathcal{F}$ has compact closure in $C(X,Y)$ with the compact open
topology.
\end{enumerate}

\begin{proof}
$1\Rightarrow 2$. See \cite{Munkres} p290. Local compactness of $X$ is not
required. However the proof makes essential use of the Tychonoff theorem,
the fact that the arbitrary product of compacta is compact.

$2\Rightarrow 1.$ Suppose in order to obtain a contradiction that $\mathcal{F%
}$ is not equicontinuous at $x.$ Choose $\varepsilon >0$ and $%
x_{n}\rightarrow x$ and $f_{n}\in \mathcal{F}$ such that $%
d(f_{n}(x),f_{n}(x_{n}))\geq \varepsilon .$ Let $A=\{x\}\cup
\{x_{1},x_{2},..\}.$ By compactness of $\overline{\mathcal{F}}$ choose a
convergent subsequence such that $f_{n_{k}}\rightarrow f.$ In particular $%
f_{n_{k}A}\rightarrow f_{A}$ uniformly and hence $f_{n_{k}A}$ is
equicontinuous. Thus $d(f_{n_{k}}(x),f_{n_{k}}(x_{n_{k}}))\rightarrow 0$ and
we have our contradiction.
\end{proof}

\subsection{Remarks}

\begin{remark}
\label{isorec}Suppose $(X,d)$ is a metric space and $f:X\rightarrow X$ is a
surjective function such that $d(x,y)=d(f(x),f(y))$ for all $x,y\in X.$ Then 
$f$ is an isometry.
\end{remark}

\begin{remark}
\label{equirec}Suppose $X$ and $Y$ are metric spaces $\mathcal{F\subset }%
C(X,Y)$ and $\mathcal{F}_{A}$ is equicontinuous for each compact subset $A.$
Then $\mathcal{F}$ is equicontinuous.
\end{remark}

\begin{remark}
\label{invariant}Suppose $h:X\rightarrow X$ is a homeomorphism, $B\subset X$
and $A=\overline{\cup _{n\in Z}h^{n}(B)}$. Then $h(A)=A.$
\end{remark}

\begin{remark}
\label{isocmpt}If $(A,d)$ is a compact metric space then $ISOM(A,d)$ is a
compact topological group.
\end{remark}

\begin{remark}
\label{metrizable}If $X$ is a locally compact separable metrizable space
then the one point compactification of $X$ is metrizable. See Theorem 8.6
p.247 of Dugundji \cite{Dug}.
\end{remark}

\begin{remark}
\label{hyygrp}Suppose $Y$ is a compact metric space. Then $H(Y,Y)$ is a
topological group under function composition in the uniform topology.
\end{remark}

\subsection{Lemmas}

\begin{lemma}
\label{orbclos}Suppose $X$ and $Y$ are metric spaces, $\mathcal{F}\subset
C(X,Y)$ is compact and $B\subset X.$ Then $A=\cup _{f\in \mathcal{F}}f(B)$
is compact.
\end{lemma}

\begin{proof}
Let $f_{n}(x_{n})\in A$ with $f_{n}\in \mathcal{F}$ and $x_{n}\in B.$ Since $%
\mathcal{F}$ is sequentially compact, and $B$ is compact. Choose convergent
subsequences $f_{n_{k}}\in \mathcal{F}$ and $x_{n_{k}}\in B$ such that $%
f_{n_{k}}\rightarrow f\in \mathcal{F}$ and $x_{n}\rightarrow x\in B.$ Hence $%
f_{n_{k}}(x_{n_{k}})\rightarrow f(x).$ Thus $A$ is sequentially compact and
hence compact.
\end{proof}

\begin{lemma}
\label{apiffeq}Suppose $(X,d)$ is a metric space, $h:X\rightarrow X$ is a
homeomorphism, and $\overline{\{h^{n}(x)\}}$ is compact $\forall x\in X.$
Then $h$ is compactly almost periodic iff $\{h^{n}\}$ is equicontinuous.
\end{lemma}

\begin{proof}
Suppose $h$ is compactly almost periodic and $B\subset X$ is compact. Let $A=%
\overline{\cup _{n\in Z}h^{n}(B)}.$ By Remark \ref{invariant} $h(A)=A.$ Thus
by Theorem \ref{Gotts} $\{h_{A}^{n}\}$ is equicontinuous. Since $B\subset A$ 
$\{h_{B}^{n}\}$ is equicontinuous. Thus

by Remark $\{h^{n}\}$ is equicontinuous. Conversely suppose $\{h^{n}\}$ is
equicontinuous and $B\subset X$ is compact. By Lemma \ref{orbclos} $A=%
\overline{\cup _{n\in Z}h^{n}(B)}$ is compact. Suppose $A\subset X$ is
compact and $h(A)=A.$ Then $\{h_{A}^{n}\}$ is equicontinuous and hence by
Theorem \ref{Gotts} $h_{A}$ is almost periodic. Thus $h$ is compactly almost
periodic.
\end{proof}

\begin{lemma}
\label{apequi}Suppose $(X,d)$ is a metric space and $h:X\rightarrow X$ is an
almost periodic homeomorphism. Then $\{h^{n}\}$ is equicontinuous.
\end{lemma}

\begin{proof}
Suppose $\varepsilon >0$ and $x\in X.$ Choose $N>0$ so that each block of $N$
iterates contains a map $h^{m}$ such that $d(id,h^{m})<\varepsilon /3.$
Choose $\delta >0$ such that $d(x,y)<\delta \Rightarrow
d(h^{i}(x),h^{i}(y))<\varepsilon $ whenever $\left| i\right| \leq N.$
Suppose $n\in Z.$ Choose $m$ such that $0\leq m\leq N$ and $%
d(h^{n-m},id)<\varepsilon /3.$ Note $d(h^{n},h^{m})<\varepsilon /3.$ If $%
d(x,y)<\delta $ then $d(h^{n}(x),h^{n}(y))\leq
d(h^{n}(x),h^{m}(x))+d(h^{m}(x),h^{m}(y))+d(h^{m}(y),h^{n}(y))$

$<\varepsilon .$ Thus $\{h^{n}\}$ is equicontinuous.
\end{proof}

\begin{lemma}
\label{tb}Suppose $(X,d)$ is a metric space and $h:X\rightarrow X$ is an
almost periodic homeomorphism. Then $\overline{\{h^{n}(x)\}}$ is totally
bounded $\forall x\in X.$
\end{lemma}

\begin{proof}
Suppose $x\in X$ and $\varepsilon >0.$ Choose $N>0$ so that each block of $N$
iterates contains a map $h^{k}$ such that $d(id,h^{k})<\varepsilon /2.$ Let $%
A=\{x,h(x),...h^{N}(x)\}.$ Suppose $y\in \overline{\{h^{n}(x)\}}.$ Let $%
d(y,h^{n}(x))<\varepsilon /2.$ Choose $m$ such that $0\leq m\leq N$ and $%
d(h^{n-m},id)<\varepsilon /2.$ Note $d(h^{n},h^{m})<\varepsilon /2.$ Thus $%
d(y,h^{m}(x))\leq d(y,h^{n}(x))+d(h^{n}(x),h^{m}(x))<\varepsilon $ and hence 
$\overline{\{h^{n}(x)\}}$ is totally bounded.
\end{proof}

\begin{lemma}
\label{equiiso}Suppose $(X,d)$ is a metric space and $h:X\rightarrow X$ is a
homeomorphism such that $\{h^{n}\}$ is equicontinuous over $(X,d)$ and $%
\overline{\{h^{n}(x)\}}$ is bounded$\forall x\in X$. Then $d^{\ast
}(x,y)=\sup_{n\in N}d(h^{n}(x),h^{n}(y))$ is an equivalent metric and $h$ is
an isometry of $(X,d^{\ast }).$
\end{lemma}

\begin{proof}
Note $d^{\ast }$ is well defined since $d$ is bounded on $\overline{%
\{h^{n}(x)\}}\cup \overline{\{h^{n}(y)\}}$

1) By definition $d^{\ast }\geq 0$ since $d\geq 0.$

2) If $x=y$ then $h^{n}(x)=h^{n}(y)\forall n\in Z$ and hence $d^{\ast
}(x,x)=0.$

3) If $x\neq y$ then $d^{\ast }(x,y)\geq d(id(x),id(y))=d(x,y)>0.$

4) $d^{\ast }$ is symmetric since $d$ is symmetric.

5) Suppose $\varepsilon >0$. Choose $n\in Z$ so that $d^{\ast
}(x,z)<d(h^{n}(x),h^{n}(z))+\varepsilon .$

Note $d(h^{n}(x),h^{n}(z))\leq d(h^{n}(x),h^{n}(y))+d(h^{n}(y),h^{n}(z))\leq
d^{\ast }(x,y)+d^{\ast }(y,z).$

Thus $d^{\ast }(x,z)\leq d^{\ast }(x,y)+d^{\ast }(y,z)+\varepsilon $ and
hence, since $\varepsilon $ is arbitrary, $d^{\ast }(x,z)\leq d^{\ast
}(x,y)+d^{\ast }(y,z).$

6)\ To show $d$ and $d^{\ast }$ are comparable it suffices to prove $%
ID:(X,d^{\ast })\rightarrow (X,d)$ is a homeomorphism. Suppose $%
x_{n}\rightarrow x.$ Then $ID(x_{n})\rightarrow ID(x)$ since $d(x,x_{n})\leq
d^{\ast }(x,x_{n}).$ Thus $ID$ is continuous. To prove $ID^{-1}$ is
continuous at $x$ suppose $\varepsilon >0$. By equicontinuity of $\{h^{n}\}$
at $x$ choose $\delta >0$ so that $d(y,x)<\delta \Rightarrow
d(h^{n}(y),h^{n}(x))<\varepsilon /2$. Suppose $d(x,y)<\delta .$ Then $%
d(h^{n}(x),h^{n}(y))<\varepsilon /2$. Thus $d^{\ast }(y,x)\leq \varepsilon
/2<\varepsilon .$ Hence $ID^{-1}$ is continuous.

7) To prove $h$ is an isometry of $(X,d^{\ast })$ it suffices by Remark \ref
{isorec} since $h$ is surjective to prove $d^{\ast }(x,y)=d^{\ast }(h^{\ast
}(x),h^{\ast }(y)).$ Note $d^{\ast }(h(x),h(y))=\sup_{n\in
Z}d(h^{n}(x),h^{n}(y))=\sup_{n\in Z}(d(h^{n+1}(x),h^{n+1}(y)))=d^{\ast
}(x,y).$

Thus $h$ is an isometry of $(X,d^{\ast }).$
\end{proof}

\subsection{Theorem \ref{Bcmptgrp}}

We assume throughout this section that $G$ is a countable abelian group of
isometries of the metric space $(X,d)$ and that $\overline{G}$ is a compact
subspace of $C(X,X).$ In general the isometry group of a metric space $X$ is
not a closed subspace of $C(X,X)$ and compact isometry groups are generally
not metrizable ( See examples \ref{nometric} and \ref{noclosed}). Hence some
care is required to establish Theorem \ref{Bcmptgrp}. The countability of $G$
is not essential for these results but enables us to detect limit points of $%
G$ with sequences rather than nets.

\begin{theorem}
\label{Bcmptgrp}Suppose $(X,d)$ is a metric space and $G\subset C(X,X)$ is a
countable abelian group of isometries of $X$ such that $\overline{G}$ is
compact in $C(X,X)$ with the compact open topology. Then $\overline{G}$ is a
compact abelian topological group of isometries of $X.$
\end{theorem}

\begin{proof}
This is a direct consequence of Lemmas \ref{gbar} and \ref{closedmult}.
\end{proof}

\begin{lemma}
\label{contmult}If $(X,d)$ is a metric space then isometry composition is
continuous in the sense that the function $\phi :ISOM(X,d)\times
ISOM(X,d)\rightarrow ISOM(X,d)$ defined as $\phi (f,g)=fg$ is continuous.
\end{lemma}

\begin{proof}
First note that $\phi $ is well defined since composition of homeomorphisms
is again a homeomorphism and $d(fg(x),fg(y))=d(g(x),g(y))=d(x,y).$ Thus if
each $\{f,g\}\subset ISOM(X,d)$ then by Remark \ref{isorec} $fg\in
ISOM(X,d). $ Suppose $U=B_{A}(fg,\varepsilon )\subset ISOM(X,d).$ Let $%
C=g(A) $. Let $V=B_{C}(f,\varepsilon /2)\times B_{A}(g,\varepsilon /2).$
Suppose $(f_{i},g_{j})\in V$ and $x\in A.$ Then $d(f_{i}g_{j}(x),fg(x))\leq
d(f_{i}g_{j}(x),f_{i}(g(x)))+d(f_{i}g(x),f(g(x))).$ Since $f_{i}$ is an
isometry, $d(f_{i}g_{j}(x),f_{i}(g(s)))=d(g_{i}(x),g(x))<\varepsilon /2.$
Since $g(x)\in C$ $d(f_{i}(g(x),fg(x)))<\varepsilon /2.$ Thus $%
d(f_{i}g_{j}(x),fg(x))<\varepsilon $ and hence $f_{i}g_{j}\in U.$ Thus $\phi 
$ is continuous at $(f,g)$ and hence continuous.
\end{proof}

\begin{lemma}
\label{orbitequals}If $A\subset X$ is compact then $\overline{G(A)}=%
\overline{G}(A)$ and $\overline{G(A)}$ is compact.
\end{lemma}

\begin{proof}
Suppose $y\in \overline{G(A)}.$ If $y\in G(A)$ then $y\in \overline{G}(A).$
Suppose $y\in \overline{G(A)}\backslash G(A).$ Let $y=\lim g_{n}(a_{n})$
where $a_{n}\in A$ and $g_{n}\in G.$ Since we may pass to convergent
subsequences we assume wolog that $a_{n}\rightarrow a\in A$ and $%
g_{n}\rightarrow g\in \overline{G}.$ Let $B=\{a\}\cup \{a_{n}\}.$ Since $B$
is compact $g_{nB}\rightarrow g_{B}$ uniformly. Thus $y=\lim
g_{n}(a_{n})=\lim g_{n}(a)=g(a).$ Hence $\overline{G(A)}\subset \overline{G}%
(A).$ Conversely, suppose $y\in \overline{G}(A).$ Let $y=g(a)$ with $a\in A$
and $g\in \overline{G}.$ Since $G$ is countable we may choose $g_{n}\in G$
such that $g_{n}\rightarrow g.$ Hence $g_{n}(a)\rightarrow g(a)=y.$ Thus $%
y\in \overline{G(A)}$. Hence $\overline{G(A)}=\overline{G}(A)$.Since $X$ is
a metric space it suffices to prove $\overline{G(A)}$ is sequentially
compact. Suppose $y_{n}\in \overline{G(A)}.$ Let $y_{n}=g_{n}(a_{n})$ where $%
g_{n}\in \overline{G}$ and $a_{n}\in A.$ By compactness of $\overline{G}$
and $A$ choose $g_{n_{i}}\rightarrow g\in \overline{G}$ and $%
a_{n_{i}}\rightarrow a\in A.$ Let $B=\{a\}\cup \{a_{n_{i}}\}.$ Since $B$ is
compact $g_{n_{i}B}\rightarrow g_{B}$ uniformly. Then $\lim y_{n_{i}}=\lim
g_{n_{i}}(a_{n_{i}})=g(a).$ Thus $\overline{G(A)}$ is compact.
\end{proof}

\begin{lemma}
\label{gbar}$\overline{G}$ is a collection of isometries of $X.$
\end{lemma}

\begin{proof}
Suppose $g\in \overline{G}$ and $y\in X.$ Since $G$ is countable choose $%
g_{n}\in G$ such that $g_{n}\rightarrow g.$ Let $y=g_{n}(a_{n})$ where $%
a_{n}\in A.$ Choose a convergent subsequence $a_{n_{i}}\rightarrow a\in A.$
Let $B=\{a\}\cup \{a_{n_{i}}\}.$ Then $y=\lim g_{n_{i}}(a_{n_{i}})=g(a).$
Thus $g$ is surjective. Suppose $x,y\in X.$ By continuity of $d$ at $%
(g(x),g(y))$ choose $\delta >0$ such that if $d(w,g(x))<\delta $ and $%
d(z,g(y))<\delta $ then $\left| d(w,g(x))-d(z,g(y))\right| <\varepsilon .$
Choose $h\in ISOM(X,d)$ so that $d(h_{\{x,y\}},g_{\{x,y\}})<\delta .$ Let $%
z=h(x)$ and $w=h(y)$ and observe that $\left| d(g(x),g(y))-d(x,y)\right|
=\left| d(g(x),g(y))-d(h(x),h(y))\right| <\varepsilon .$ Thus $%
d(x,y)=d(g(x),g(y)).$ Letting $Y=im(g)$ and $\delta =\varepsilon $ shows
that $g$ is a homeomorphism of $X$ onto $Y.$
\end{proof}

\begin{lemma}
\label{closedmult}Suppose $f,g\in \overline{G}$. Then $fg\in \overline{G},$ $%
g^{-1}\in \overline{G},$ $\overline{G}$ is abelian, and the map $\phi :%
\overline{G}\rightarrow \overline{G}$ is continuous where $\phi (g)=g^{-1}.$
\end{lemma}

\begin{proof}
Since $G$ is countable let $f_{n}\rightarrow f$ and $g_{n}\rightarrow g$
with $f_{n},g_{n}\in G.$ Then $f_{n}g_{n}\rightarrow fg$ by Lemma \ref
{contmult} and hence $fg\in \overline{G}$ since $f_{n}g_{n}\in G$ and $%
\overline{G}$ is closed. By Lemma \ref{gbar} $g$ is an isometry and hence $%
g^{-1}$ exists. We will prove $g_{n}^{-1}\rightarrow g^{-1}$ and conclude
that $g^{-1}\in \overline{G}$ since $\overline{G}$ is closed. Suppose $%
B\subset X$ is compact. Let $A=\overline{G}(B).$ By Lemma \ref{orbitequals} $%
A$ is compact and invariant under $G.$ Hence $g_{n|A}\rightarrow g_{|A}$ and
thus by Remark \ref{isocmpt} $(g_{n|A})^{-1}\rightarrow
(g_{|A})^{-1}=(g^{-1})_{|A}.$ Thus $g^{-1}\in \overline{G}.$ Hence $\phi $
is well defined. To check continuity let $U=B_{B}(g^{-1},\varepsilon )\cap 
\overline{G}\subset \overline{G}.$ Let $A=\overline{G}(B).$ Let $%
V=B_{A}(g,\varepsilon )\cap \overline{G}.$ Note $\phi (V)\subset U.$ Thus $%
\phi $ is continuous at $g$ and hence continuous. By Lemma \ref{contmult} $%
fg(x)=\lim f_{n}g_{n}(x)=\lim g_{n}f_{n}(x)=gf(x).$ Thus $\overline{G}$ is
abelian.
\end{proof}

\end{document}